\documentclass [a4paper,11pt]{article}

\usepackage{amsthm}
\usepackage{amsmath,amssymb,latexsym,amsfonts,mathrsfs}
\usepackage[dvipdfmx]{graphicx}

\usepackage{color}

\usepackage{mathtools}
\usepackage{fancyhdr}
\usepackage[top=30truemm, bottom=30truemm, left=25truemm, right=25truemm]{geometry}

\newcommand{\ep}{\varepsilon}
\newcommand{\nn}{\nonumber}

\newcommand{\SCR}[1]{{\mathscr #1}}

\newcommand{\CAL}[1]{{\cal #1}
}

\newcommand{\J}[1]{\left\langle #1 \right\rangle}
\newcommand{\D}[1]{{\mathscr D}( #1 )}

\theoremstyle{definition}
\newtheorem{Thm}{{\bf Theorem}}[section]

\newtheorem{Lem}[Thm]{{\bf Lemma}}
\newtheorem{Prop}[Thm]{{\bf Proposition}}

\newtheorem{Ass}[Thm]{{\bf Assumption}}

\newtheorem{Rem}[Thm]{{\bf Remark}}

\newcounter{Exami}

\newcommand{\Proof}[2][Proof]{
\begin{proof}[{\bf #1}]
#2
\end{proof}
}

\begin{document}

\begin{flushleft}
{\bf \Large Final state problem for nonlinear Schr\"{o}dinger equations with time-decaying harmonic oscillators
} \\ \vspace{0.3cm} 
by {\bf \large Masaki Kawamoto } \\
Department of Engineering for Production, Graduate School of Science and Engineering, Ehime University, 3 Bunkyo-cho Matsuyama, Ehime 790-8577. Japan \\
Email: {kawamoto.masaki.zs@ehime-u.ac.jp} 
\end{flushleft}

\begin{center}
\begin{minipage}[c]{400pt}
{\bf Abstract}. {\small
We consider the final-state problem for the nonlinear Schr\"{o}dinger equations (NLS) with a suitable time-decaying harmonic oscillator. In this equation, the power of nonlinearity $|u|^{\rho}u $ is included in the long-range class if $0 < \rho \leq  2/(n(1- \lambda)) $ with $0 \leq \lambda <1/2$, which is determined by the harmonic potential and a coefficient of Laplacian. In this paper, we find the final state for this system and obtain the decay estimate for asymptotics.  
}
\end{minipage}
\end{center}

\begin{flushleft}
{\bf MSC 2010}; Primary: 35Q55, Secondly: 35J10
\end{flushleft}

\begin{flushleft}
{\bf Keywords}; Nonlinear scattering theory; Long-range scattering; Time-dependent harmonic oscillators; Final state problem; 
\end{flushleft}
\section{Introduction}
{In this paper}, we consider the nonlinear Schr\"{o}dinger equations with time-dependent harmonic potentials. 
\begin{align}\label{eq1}
\begin{cases}
i \partial _t u(t,x) - \left( - \Delta /2 + \sigma (t) |x|^2/2 \right) u(t,x) = \mu |u(t,x)|^{2/(n(1- \lambda))} u(t,x) , \\
u(0,x) = u_0 (x),
\end{cases}
\end{align}
where $(t,x) \in {\bf R} \times {\bf R}^n$, $n \in \{ 1,2,3\}$, and $ \mu  \in {\bf R} \backslash \{0 \}$. $\lambda \in [0,1/2)$ is defined later. We let
\begin{align*}
H_0 (t) := - \Delta /2 + \sigma (t) |x|^2/2 
\end{align*} 
and employ the following assumption on the coefficient of the harmonic oscillator $\sigma (t)$. 
\begin{Ass} \label{A1}
Let $\zeta _1 (t)$ and $\zeta _2 (t)$ be the solutions to 
\begin{align*}
\zeta _j ''(t) + \sigma (t) \zeta _j (t) =0, \quad
\begin{cases}
\zeta _1 (0) = 1, \\
\zeta _1 ' (0) =0,
\end{cases}
\quad
\begin{cases}
\zeta _2 (0) = 0, \\
\zeta _2 ' (0) =1.
\end{cases}
\end{align*} 
Then, there exist $r_0 >0$, $c_{1, \pm} , c_{2, \pm} \notin \{ 0, \infty , - \infty\}$, $c_{3,\pm} \in {\bf R}$ and $\lambda \in [0, 1/2 )$ such that for all $|t| > r_0$, the followings hold: 
$$
|\zeta _2 (t)| \geq c
$$
and 
\begin{align*}
\lim_{t \to \pm \infty} \frac{\zeta _1 (t)}{(\pm t)^{ \lambda} } = c_{1, \pm}, \quad 
\lim_{t \to \pm \infty} \frac{\zeta _2 (t)}{(\pm t)^{1-\lambda}} = c_{2, \pm}, \quad 
\lim_{t \to \pm \infty} \frac{\left|\zeta _2 (t) - c_{2, \pm}(\pm t)^{1-\lambda}\right|}{(\pm t)^{\lambda}} = c_{3, \pm}.
\end{align*}
Moreover, $\zeta_1 (t)$, $\zeta_2 (t)$, $\zeta_1' (t)$, and $\zeta_2' (t)$ are continuous functions.
\end{Ass}
This assumption implies $\zeta _1 (t)  \sim c_{1, \pm}|t|^{\lambda} + o(|t|^{\lambda})$ and $\zeta _2 (t) = c_{2, \pm} |t|^{1- \lambda} + \CAL{O}(|t|^{ \lambda}) $. If $\sigma (t)$ decays in $t$, there are some examples satisfying this assumption, see, e.g., Kawamoto-Yoneyama \cite{KY} and references their in. Based on this assumption, the quantum particle governed by the energy $H_0(t)$ is decelerated with velocity $v (t) = \CAL{O} (|t|^{- \lambda})$ as $|t| \to \infty $; however, it is not trapped (see Kawamoto \cite{Ka}, Ishida-Kawamoto \cite{IK}). This phenomenon changes the threshold between the long range and the short range of the power of nonlinearity from $2/n$ to $2/n(1- \lambda)$. This was discovered by Kawamoto-Muramatsu \cite{KM}. In this study, we investigate the nature of the final states upon inclusion of nonlinearity in the long-range class. 

We first consider $U_0(t,s)$ to be a propagator for $H_0 (t)$, that is, the family of unitary operators $\{ U_0(t,s)\}_{(t,s) \in {\bf R}^2} $ acts on $L^2({\bf R}^n)$ with conditions that for all $t,s,\tau \in {\bf R}$, the following hold for $\D{H_0(s)}$.
\begin{align*}
& i \partial _t U_0(t,s) = H_0(t) U_0(t,s), \quad i \partial _s U_0(t,s) = -U_0(t,s) H_0(s), \\
& U_0(t, \tau) U_0(\tau, s) = U_0(t,s), \quad U_0(s,s) = \mathrm{Id}_{L^2({\bf R}^n)}, \quad U_0(t,s) \D{H_0(s)} \subset \D{H_0(s) }. 
\end{align*}
Using $\zeta _1 (t)$ and $\zeta _2 (t)$, we obtain the following \CAL{MDFM-decomposition}, see Korotyaev \cite{Ko} and Carles \cite{Ca} (also see Kawamoto \cite{Ka2} and \cite{KM}).

\begin{Lem}
For $\phi \in \SCR{S}({\bf R}^n)$, we define
\begin{align*}
\left( \CAL{M}(\tau) \phi \right) (x) = e^{ix^2/(2 \tau)} \phi (x), \quad
\left(
\CAL{D}(\tau) \phi
\right) (x) = \frac{1}{(i \tau)^{n/2}} \phi (x/ \tau).
\end{align*}
Then, the following MDFM decomposition holds. 
\begin{align}\label{mdfm}
U_0(t,0) = \CAL{M} \left(  \frac{\zeta _2(t)}{\zeta _2 '(t)} \right) \CAL{D} (\zeta _2 (t)) \SCR{F} \CAL{M} \left( \frac{\zeta _2 (t)}{ \zeta _1 (t)} \right) , 
\end{align}
where $\SCR{F}$ indicates the Fourier transform.
\end{Lem}
\begin{Rem}
If $\sigma (t) = \sigma _0 t^{-2}$ for some $\sigma _0 \in (0,1/4)$ and for all $|t| \geq r_0$, then $y_1(t)$ and $y_2(t)$, the linearly independent solutions of $y''(t) + \sigma (t) y(t) =0$, satisfy $y_1 (t) = |t|^{\lambda}$ and $c_2 (t) = |t|^{1- \lambda}$ for $|t| \geq r_0$ with $\lambda = (1 -\sqrt{1-4 \sigma _0})/2$. In order to gain the time decay, we decompose $$ \CAL{M} \left( {\zeta _2 (t)}/{ \zeta _1 (t)} \right) = \CAL{M} \left( {\zeta _2 (t)}/{ \zeta _1 (t)} \right) - 1 + 1 $$ 
and use $ | \CAL{M} \left( {\zeta _2 (t)}/{ \zeta _1 (t)} \right) - 1 | \leq C |\zeta _2(t)/\zeta _1 (t)|^{-\theta} |x|^{2 \theta} $, $0 \leq \theta \leq 1$, which acts critical role in analyzing long-range nonlinearity. Under the assumption \ref{A1}, one can get the time decay such that $|\zeta _2(t)/\zeta _1 (t)|^{-\theta} \leq C t^{- \theta (1- 2 \lambda)}$. However, as usual there is a possibility to be that $\zeta _1 (t) = \CAL{O} (t^{1- \lambda})$ or $\zeta _2 (t) = \CAL{O} (t^{\lambda})$. In such case we can not get the time decay from $\CAL{M} \left( {\zeta _2 (t)}/{ \zeta _1 (t)} \right) - 1$ and our scheme fails. In this sense, we need to assume Assumption \ref{A1}.
\end{Rem}
Hereafter, we use the notation 
\begin{align*}
\CAL{M}_1 (t) = \CAL{M} \left(  \frac{\zeta _2(t)}{\zeta _2 '(t)} \right), \quad \CAL{M}_2 (t) =\CAL{M} \left( \frac{\zeta _2 (t)}{ \zeta _1 (t)} \right),
\end{align*}
$u(t, \cdot ) = u(t)$, and 
\begin{align*}
F(u(t)) = \mu |u(t)|^{2/n(1- \lambda)} u(t).
\end{align*}
We let 
\begin{align*}
u_p (t) &:= \CAL{M}_1(t) \CAL{D}(\zeta _2 (t)) \hat{w} (t), \\ 
\hat{w} (t)  &:= \widehat{u_+} \mathrm{exp} \left(  {-i \mu |\widehat{u_+}| ^{2/n(1- \lambda)} ( \log t) /c_+} \right), 
\end{align*}
where $\hat{\cdot}$ denotes the Fourier transform and $c_+ := |c_{2,+}|^{1/(1- \lambda)}$. Because of some technical reasons, the following restriction on $\lambda $ is mandatory. We suppose that $\lambda = \lambda (n)$ satisfies 
\begin{align}\label{15}
0 \leq \lambda (1) <  \frac{13-\sqrt{145}}{4} \quad 0 \leq \lambda (2) < \frac{7-\sqrt{41}}{4} , \quad 0 \leq \lambda (3) < 0.022. 
\end{align}
For simplicity, we set constants $k$ and $\alpha$ as 
\begin{align*}
0< k  <  1 + \frac{2}{n(1- \lambda)}, \quad 0 < \alpha < \mathrm{min} ( 1, k/2 -n\lambda /4 ), 
\end{align*}
and $H^{\gamma, \nu} = \{ f \in \SCR{S}'({\bf R}^n) \, ; \, \| f \|_{\gamma, \nu} = \left\| (1+ |x|^2)^{\nu /2} (1- \Delta)^{\gamma /2} f \right\|_{L^2({\bf R}^n)} < \infty \}$. Then, we have the following theorem: 
\begin{Thm}\label{T}
Under the assumption \ref{A1}, we suppose that $\lambda = \lambda (n)$ satisfies \eqref{15}. Further, we suppose that $\widehat{u_+} \in H^{k,0}\cap L^{\infty} ({\bf R}^n)$ with $\| \widehat{u_+}   \|_{\infty} \ll 1 $. Then, there exists $T>0$ such that a solution $u(t)$ in \eqref{eq1} exists globally in $C ([T, \infty) \, ; \, L^2({\bf R}^n))$ and that for all 
\begin{align}\label{16}
 \frac{n(-2 \lambda ^2 + \lambda +1) + 8 \lambda}{4}  < b < \lambda + \alpha (1- 2 \lambda) , 
\end{align}
the inequality
\begin{align*}
\sup_{t \in [T, \infty )} t^b \left\| u(t) - u_p (t)  \right\|_{L^2({\bf R}^n)} < \infty 
\end{align*}
holds. 
\end{Thm}
If $\sigma (t) \equiv 0$, \eqref{eq1} corresponds to the nonlinear Schr\"{o}dinger equations with long-range nonlinearity. The asymptotic behavior and the final state problem have been considered in several studies, e.g., Ginibre-Ozawa \cite{GO}, Hayashi-Ozawa \cite{HO}, Ozawa \cite{O}, Hayashi-Naumkin \cite{HN}, Hayashi-Naumkin-Wang \cite{HNW}, Masaki-Miyazaki \cite{MM}, Masaki-Miyazaki-Uriya \cite{MMU}. In the case of $\lambda = 0$, which includes the case where $\sigma (t) \equiv 0$, the condition for $b$ corresponds to (or includes) that of \cite{MMU} and \cite{MM}. Hence, our result will be a natural extension to the final state problem for long-range NLSs with harmonic potentials. On the other hand, in the previous studies, nonlinearities have been admitted to have some {\em non-resonant part} and we have not dealt with such nonlinearities with our model. In order to deal with such terms, a strong restriction for $\alpha _n$ in the auxiliary space $L^{\beta_n }([T, \infty);L^{\alpha _n} ({\bf R}^n))$ (see \S{2.2}) may be necessary, and hence a more rigorous argument for $\lambda$ may appear. For simplicity, we avoided these issues in this paper; however, our result will be the first step toward considering these issues.

By the result in \cite{HN} and \cite{KM}, one can have for all $\lambda \in (0,1/2)$ and some $u_0 \in H^{0, \gamma} \cap H^{\gamma ,0}$ with $\| u_0\|_{\gamma,0} + \| u_0 \|_{0, \gamma} < \ep \ll 1$, there exists a function $W$ and $\Psi  $ such that 
\begin{align*}
\left\|  u (t) -e^{-i \mu |W|^{2/(n(1- \lambda))} \log t + i \Psi } U_0(t,0) \SCR{F}^{-1} W
\right\|_{k} \leq C t^{-(1-2\lambda) \alpha_0 + C(\ep)}, 
\end{align*}
with $\alpha_0 < \min (1/2 + 1/(n(1- \lambda)) -n/4, 1) $, where $k= 2$ or $\infty$. In this result, one can deal with all $\lambda \in (0,1/2)$ but the decay $-(1-2\lambda) \alpha_0 $ is weak compared with $b$. Hence our result is a more rigorous analysis for asymptotic behavior of $u(t)$. 
We use the approach established by \cite{HNW}. However, in order to imitate this approach, a strong restriction for $\lambda$ is necessary. As an ideal case, we consider $\lambda$ arbitrarily in $[0,1/2)$. Hence, in this sense, there is significant scope for development.

In section \S{2}, we first decompose $u(t) - u_p (t)$ through the MDFM-decomposition and deduce the similar formula according to \cite{HNW};
\begin{align}\label{ad-1}
u(t) - u_p(t) = i \int_t^{\infty} U_0(t,s) \left( 
F(u(s)) - F(u_p(s))
\right)  ds + \CAL{E} (t) + \CAL{A}(t),
\end{align} 
where $\CAL{A}(t)$ is the special term appearing only this model but to estimate this term is easy (see, Lemma \ref{L-ad1}). Next we consider the fractional derivative for $\hat{w}(t)$. Since $2/(n(1-\lambda)) $ is not in ${\bf N}$ even if $n=1,2$, we have to consider the fractional derivative for any dimension. At last, we set energy spaces in order to construct the contraction mapping. In \S{3}, we prove Theorem \ref{T}. In analyzing the first term of r.h.s of \eqref{ad-1}, the lower restriction for $b$ appears. On the other hand, in analyzing $\CAL{E}(t)$ and $\CAL{K}(t)$, the upper restriction for $b$ appears. In order to justify the lower and upper restriction for $b$, the restriction for $\lambda$ is demanded, see \S{3.2}. 


\section{Preliminaries}
Throughout, for $q \in [1,\infty]$, the norm $\| \cdot \|_{q}$ denotes the norm on $L^q({\bf R}^n)$. Additionally, for $q,r \in [1,\infty]$ and $a<b$, we define the time-weighted Bochner-Lebesgue space $L^q_{\lambda} ((a,b) ; L^r({\bf R}^n))$ (see Kawamoto-Yonenayama \cite{KY} and Kawamoto \cite{Ka3}) as follows: 
\begin{align*}
L_{\lambda}^q ((a,b) ; L^r ({\bf R}^n)) = \left\{ 
F \in \SCR{S}' ((a,b) \times {\bf R}^{n}) \, \middle| \, 
 \| F \|_{L^q ((a,b) ; L^r({\bf R}^n), \lambda )} < \infty
\right\},
\end{align*}
where 
\begin{align*}
\| F \|_{L^q ((a,b) ; L^r({\bf R}^n) \lambda )} := \left( 
\int_{a}^{b} (1 + t^2)^{- \lambda /2} \left\| 
F(t, \cdot)
\right\|_{r}^q dt 
\right) ^{1/q} .
\end{align*}

We say that the pair $(q,r)$ is admissible if it satisfies 
\begin{align*}
\frac{1}{q} + \frac{n}{2r} = \frac{n}{4}, \quad q >2 , \quad r \geq 2.
\end{align*}
Then, we introduce the so-called Strichartz estimate for associated propagator $U_0(t,s)$, which was first obtained by \cite{KY} for the case of the limited coefficients $\sigma (t)$ and next under more generalized condition, including assumption \ref{A1}, which was considered by \cite{Ka3}. 
\begin{Lem}
Let $(q,r)$ and $(\tilde{q}, \tilde{r} )$ be admissible pairs, and let $ s' $ denote the H\"{o}lder exponent of $s$, i.e., $1/s + 1/s' =1$.
Then, for $\phi \in L^2({\bf R}^n)$ and $F \in L_{\lambda} ^{\tilde{q}'}((a,b) ; L^{\tilde{r}'}({\bf R}^n), - \lambda)$, there exists $C>0$ that is independent of $a,b$ such that the following inequalities hold.
\begin{align*}
\left\| 
U_0(t,0) \phi
\right\|_{L^q((a,b);L^r({\bf R}^n) , \lambda)} \leq C \| \phi \|_2
\end{align*}
and 
\begin{align*}
\left\| 
\int_{s<t} U_0(t,s) F(s) ds 
\right\|_{L^q ((a,b);L^r({\bf R}^n), \lambda)} \leq C \left\| F \right\|_{L^{\tilde{q}'}((a,b);L^{\tilde{r}'}({\bf R}^n) , - \lambda)}
\end{align*}
\end{Lem}
In the following, we also assumed that $\lambda >0$.
\subsection{Setting}
In this study, we use the decomposition formula given by Hayashi-Naumkin-Wang \cite{HNW} (also see \cite{MM} and \cite{MMU}). Throughout we assume that $c_{2, +} >0$, i.e., $\zeta _2 (s) >0$, $s \geq r_0$ for simplicity. The case where $c_{2,+} <0$ can be handled by the same way. We let
\begin{align*}
\CAL{L} (f(t)) = (i \partial _t + \Delta /2 - \sigma (t) |x|^2 /2 ) f(t). 
\end{align*}
Then, by the definition of $u(t)$, we obtain 
\begin{align*}
\CAL{L} (u(t) - u_p(t)) = F(u(t)) -F(u_p(t)) - \CAL{L} (u_p(t)) + F(u_p(t)).
\end{align*}
 Using Duhamel's formula, 
\begin{align} \label{1}
u(t) - u_p(t) = i \int_t^{\infty} U_0(t,s) \left( 
F(u(s)) - F(u_p(s))
\right)  ds + \CAL{E} (t) + \CAL{A}(t), 
\end{align}
holds where we define 
\begin{align*}
\CAL{A} (t) &:=  i\int_t^{\infty} 
 U_0(t,0) \SCR{F}^{-1} \left(  \frac{c_{+}s}{\zeta _2 (s)^{1/(1- \lambda)} }-1 \right) F(\hat{w}(s)) \frac{ds}{c_+ s}, \\
\CAL{E} (t) &:= R(t) \hat{w}(t) - i \int_t^{\infty} U_0(t,s) R(s) F(\hat{w}(s)) \frac{ds}{\zeta _2 (s)^{1/(1- \lambda)}}  
\end{align*}
and
\begin{align*}
R(t) := \CAL{M}_1 (t) \CAL{D}(\zeta _2 (t)) \left(  \SCR{F}\CAL{M}_2 (t) \SCR{F}^{-1}  -1 \right).
\end{align*}
Because the following equality holds
\begin{align*}
& i \int_t^{\infty} U_0(t,s ) \left( 
- \CAL{L}(u_p(s)) + F(u_p (s)) 
\right) ds \\ & = 
i \int_t^{\infty} \Big( 
-U_0(t,0) i \partial _s U_0(0,s) \CAL{M}_1 (s) \CAL{D} (\zeta _2 (s)) \hat{w}(s)   \\  & \quad \qquad + \frac{1}{\zeta _2(s)^{1/(1-\lambda)}} U_0(t,s)\CAL{M}_1 (s) \CAL{D}(\zeta _2 (s)) F (\hat{w}(s))
\Big)  ds \\ &=  
i \int_t^{\infty} \left( 
-U_0(t,0) i \partial _s U_0(0,s) \CAL{M}_1 (s) \CAL{D} (\zeta _2 (s)) \hat{w}(s) +  U_0(t,s) U_0(s,0) \SCR{F}^{-1} \frac{F (\hat{w}(s))}{\zeta _2 (s)^{1/(1- \lambda)} } 
\right)  ds  \\ & \quad + 
i \int_t^{\infty} U_0(t,s)  \CAL{M}_1 (s) \CAL{D}(\zeta _2 (s)) \left( 1-\SCR{F}\CAL{M}_2(s) \SCR{F}^{-1} \right)  \frac{F (\hat{w}(s))}{\zeta _2 (s)^{1/(1- \lambda)} }
 ds \\ & \equiv 
 I_1 (t) + I_2(t). 
\end{align*}
We note that 
\begin{align*}
\frac{F(\hat{
w} (s)) }{c_{+} s } = i\frac{d}{ds} \hat{w} (s)
\end{align*}
with $c_+ = c_{2,+}^{1/(1- \lambda)}$ and hence get 
\begin{align*}
& i \int_t^{\infty} 
 U_0(t,s) U_0(s,0) \SCR{F}^{-1} \frac{F (\hat{w}(s))}{\zeta _2 (s)^{1/(1- \lambda)} }  ds \\ &= - \int_t^{\infty} 
 U_0(t,0) \SCR{F}^{-1} \frac{c_{+}s}{\zeta _2 (s)^{1/(1- \lambda)} } \frac{d}{ds} \hat{w} (s)  ds 
 \\ &= - \int_t^{\infty} 
 U_0(t,0) \SCR{F}^{-1}\frac{d}{ds} \hat{w} (s)  ds - \int_t^{\infty} 
 U_0(t,0) \SCR{F}^{-1} \left(\frac{c_{+}s}{\zeta _2 (s)^{1/(1- \lambda)} } -1 \right)\frac{d}{ds} \hat{w} (s)  ds. 
\end{align*}
We have 
\begin{align*}
I_1(t) = & \lim_{R \to \infty}  U_0(t,0) \left( 
U_0(0,R) \CAL{M}_1 (R) \CAL{D}(\zeta _2 (R)) \hat{w} (R) - \SCR{F}^{-1} \hat{w} (R) 
\right)  \\ & + \CAL{M}_1(t) \CAL{D}(\zeta _1(t)) \left( \SCR{F} \CAL{M}_2 (t) \SCR{F}^{-1} -1 \right) \hat{w}(t) + \CAL{A}(t)
\end{align*}
and 
\begin{align*}
I_2(t) =  
i \int_t^{\infty} U_0(t,s)  \CAL{M}_1 (s) \CAL{D}(\zeta _2 (s)) \left( 1-\SCR{F}\CAL{M}_2(s) \SCR{F}^{-1} \right)  \frac{F (\hat{w}(s))}{\zeta _2 (s)^{1/(1- \lambda)}}
 ds .
\end{align*}
\begin{Lem} \label{L1}
Let $0 < \gamma <  1+ 2/(n(1- \lambda)) $. There are constants $C >0 $, $\kappa _{1,n}(\gamma) \geq 0$, and $ \kappa _{2,n}(\gamma) \geq 1$ such that for $ |s| \geq r_0$, 
\begin{align*}
\left\| 
\hat{w} (s)
\right\| _{\gamma ,0} \leq C (\log s )^{[\gamma ]} \J{\| \widehat{u_+} \|_{{\gamma , 0}}} ^{\kappa _{1,n} (\gamma)} \J{ \| \widehat{u_+} \|_{\infty}}^{\kappa _{2,n} (\gamma)} \left\| \widehat{u_+} \right\|_{\gamma ,0}
\end{align*}
holds, where $[\gamma]$ indicates the smallest integer $n$ such that $n \geq \gamma$.
\end{Lem}
\Proof{
We imitate the approaches proposed by \cite{MMU} (see the proof of Lemma 2.6.). We first consider the case where $n=3$. By the condition $\gamma (3) < 0.022$, it is sufficient to consider the case $1+ 2/(3(1- \lambda)) < 2$. Let $\tau := \gamma -1 < \tau_0 =: 2/(3(1- \lambda)) $, $\J{\cdot} := (1 + \cdot ^2)^{1/2}$, $\widehat{u_+} = U$, and $e^{-i \mu |U|^{\tau_0} (\log t)/c_+} = \phi_U$. Then 
\begin{align*}
\left\| \J{\nabla}^{1+ \tau} \phi_U U  \right\|_2 \leq C \sum_{j=1}^3 \left\| \J{\partial_j}^{1 + \tau} \phi_U U  \right\|_2
 \leq C \sum_{j=1}^3 \left\| \J{\partial_j}^{\tau} (1 + \partial _j) \phi_U U  \right\|_2, 
\end{align*}
where we use $\| \J{\partial _j} (1 + \partial _j)^{-1} \|_{\SCR{B}(L^2 \to L^2)} \leq C $. By the tentative calculation, we have 
\begin{align*}
\partial _j \phi_U U \sim ( \log t)  \left( (\partial _j \bar{U} ) U^2 + (\partial _j U) |U|^2 \right) |U|^{\tau_0-2} \phi_U + (\partial _j U) \phi_U. 
\end{align*}
Hence, we estimate the two terms as: 
\begin{align*}
J_1 := ( \log t) \left\| D_j^{\tau} U_j |U|^{\tau_0}  \phi_U   \right\|_2 , \quad J_2:= \left\| D_j^{\tau} U_j \phi_U \right\|_2,
\end{align*}
where $D_j :=(- \partial _j ^2)^{1/2}$ and $U_j = (\partial _j U)$. By the Kato-Ponce inequality and the proposition A.1 of Visan \cite{V}, for $\tau < s <1$ and $1/p+1/q =1/2$, $p,q >2$, we have, 
\begin{align*}
J_1 &\leq C ( \log t)  \left( \left\| D_j^{\tau} |U|^{\tau_0} U_j  \right\|_2 + \left\| |U|^{\tau_0} U_j \right\|_{p} \left\| D_j^{\tau} \phi_U\right\|_{q} \right) \\ & \leq 
C ( \log t)  \left( \left\| U \right\|_{\infty}^{\tau_0} \left\| U \right\|_{\gamma, 0} + (\log t) \left\| U \right\|_{\infty} ^{\tau_0} \left\| U_j \right\|_{p} \times\left\| U \right\|_{\infty}^{\tau_0-\tau/s} \left\| D_j^s  U \right\|_{q\tau/s}^{\tau/s} \right) .
\end{align*}
Here, we also use interpolation Gagliardo-Nirenberg inequality and obtain for 
\begin{align*}
\frac{1}{p}= \frac13 + a_1 \left( \frac12 - \frac{1+ \tau}{3} \right), \quad  \frac{s}{q \tau} = \frac{s}{3} + a_2 \left( \frac12 - \frac{1+ \tau}{3} \right)  ,
\end{align*}
the inequality
\begin{align*}
J_1 &\leq C ( \log t)  \left( \left\| U \right\|_{\infty}^{\tau_0} \left\| U \right\|_{\gamma, 0} + (\log t) \left\| U \right\|_{\infty} ^{2\tau_0 - \tau /s + (2-a_1 - \tau a_2 /s)  } \left\| U \right\|_{\gamma ,0}^{a_1 + \tau a_2 /s} \right),  \\ 
& \leq C ( \log t)^2 \J{\left\| U \right\|_{\infty}}^{2\tau_0 + (1-\tau)}  \left\| U \right\|_{\gamma ,0}. 
\end{align*}
where we assume $a_1 + \tau a_2 /s=1$ and use $-\tau /s < -\tau $. Next, we consider the case where $n=1$. The case of $n=2$ can be handled in a similar manner; hence, we omit the proof for $n=2$. Assume that $\gamma >3$. We set $\tau = \gamma -3$ and $\tau _1 = 2/(1- \lambda) >2$. Similar to the case of $n=3$, we first calculate $\partial _j ^3 \phi_U U$, and notice that it is enough to deal with the most effective term
\begin{align*}
J_3 := (\mu \log t)^3 \left\| D_j^{\tau} |U|^{3 \tau_1}U_j^3 \phi_U \right\|_2
\end{align*}
since other terms only yield terms smaller than $O((\log t)^3)$. By the Kato-Ponce inequality and the inequality of \cite{V}, the following inequalities hold.
\begin{align*}
J_3 &\leq C ( \log t)^3 \left( \left\| D_j^{\tau} |U|^{3 \tau _1} U_j^3  \right\|_2 + \left\| |U|^{3 \tau _1} U_j^3  \right\|_p \left\| D_j^{\tau_1} \phi_U\right\|_q \right) 
\\ & \leq 
 C ( \log t)^3 \left[ \left\| D_j^{\tau} |U|^{3 \tau_1} \right\|_2 \| U_j \|_{\infty}^3 + \| {U} \|_{\infty}^{3 \tau_1} \left\| D_j^{\tau} U_j^3 \right\|_2   \right.
  \\   & \quad \qquad \qquad
   \left.+ (\log t)\left(  \left\| U \right\|_{\infty} ^{3\tau_1} \left\| U_j^3 \right\|_{p} \left\| U \right\|_{\infty}^{\tau_1-\tau/s} \left\| D_j^s  U \right\|_{q\tau/s}^{\tau/s}  \right)\right]. 
 \\ &  \leq 
 C ( \log t)^3 \left[ \left\| U \right\|_{\infty}^{3 \tau_1 +2 -6/(5+2 \tau)} \left\| U \right\|_{\gamma ,0}^{1 + 6/(5+ 2 \tau)} + \| U \|_{\infty}^{3\tau_1 +2-2/(5+2\tau)} \left\| U \right\|_{\gamma ,0}^{1 + 2/(5+ 2 \tau)} \right. 
 \\   & \quad  \quad \qquad \qquad
  \left. + ( \log t)\left(  \left\| U \right\|_{\infty} ^{4\tau_1 -\tau /s + 3}   \left\| U \right\|_{\gamma ,0}  \right)\right]. 
\end{align*}

}

In the following, we set 
\begin{align*}
\CAL{E}_{\gamma} (\widehat{u_+}) := \J{ \| \widehat{u_+} \|_{{\gamma , 0}} }^{\kappa _{1,n}(\gamma)} \J{ \| \widehat{u_+} \|_{\infty}}^{\kappa _{2,n}(\gamma)}  \| \widehat{u_+} \|_{{\gamma , 0}}.
\end{align*}

\begin{Lem}
Suppose that $ \widehat{ u _+}  \in H^{\gamma ,0} $ holds. Then, for all $s \in {\bf R} $ with $|s| \geq r_0$ and $0< \gamma_0 \leq  \mathrm{min} (2, \gamma) $, 
\begin{align*}
\left\| 
\left( 
\SCR{F} \CAL{M} _2 (s) \SCR{F} ^{-1} -1
\right) \hat{w} (s)
\right\|_2 \leq C \left| \frac{ \zeta _1 (s)}{\zeta _2 (s)} \right|^{\gamma_0 /2} (\log s)^{[\gamma_0]} \CAL{E}_{\gamma_0} (\widehat{u_+})
\end{align*} 
holds.
\end{Lem}
\Proof{ 
Noting that for all $\phi \in \D{|x|^{\gamma_0}}$, $ | (\CAL{M}_2 (s) -1) \phi | \leq  C | \zeta _1 (s) |^{\gamma_0 /2} |\zeta _ 2 (s)|^{-\gamma _0 /2}  | |x|^{\gamma_0} \phi |   $ holds, we have
\begin{align*}
& \left\| \left(  \SCR{F} \CAL{M}_2 (t) \SCR{F}^{-1} -1 \right) \hat{w}(t) \right\| _2
= 
 \left\| \left(  \CAL{M}_2 (t) -1 \right) \SCR{F}^{-1} \hat{w}(t) \right\|_2 \\ & \leq 
C \left| \frac{\zeta _1 (t)}{\zeta _2 (t)} \right| ^{\gamma_0/2} \left\| 
 \SCR{F}^{-1} \hat{w}(t) 
 \right\| _{{0,\gamma_0}}  \leq 
C \left| \frac{\zeta _1 (t)}{\zeta _2 (t)} \right| ^{\gamma_0/2} \left\|
  \hat{w}(t) 
 \right\| _{{\gamma_0},0} 
 \\ &  \leq 
C \left| \frac{\zeta _1 (t)}{\zeta _2 (t)} \right| ^{\gamma_0/2} (\log t)^{[\gamma _0]}
\CAL{E}_{\gamma_0} (\widehat{u_+})
,  
\end{align*}
 where we employ lemma \ref{L1}.

 }
 
The first term of $I_1(t)$ can be written as 
\begin{align*}
& \lim_{R \to \infty}  U_0(t,0) \left( 
U_0(0,R) \CAL{M}_1 (R) \CAL{D}(\zeta _2 (R)) \hat{w} (R) - \SCR{F}^{-1} \hat{w} (R)
\right) \\ & = 
 \lim_{R \to \infty}  U_0(t,0) \left( 
\CAL{M}_2 (R) ^{-1 }  - 1\right) \SCR{F}^{-1} \hat{w} (R),  
\end{align*}
and together with 
\begin{align*}
& \left\| 
 U_0(t,0) \left( 
U_0(0,R) \CAL{M}_1 (R) \CAL{D}(\zeta _2 (R)) \hat{w} (R) - \SCR{F}^{-1} \hat{w} (R) 
\right) 
\right\|_2 \\ &
 \leq C  \left| \frac{\zeta _1 (R)}{\zeta _2 (R)} \right| ^{\gamma_0/2} (\log R)^{[\gamma _0]} \CAL{E}_{\gamma _0} (\widehat{u_+})\to 0 
\end{align*}
as $R \to \infty$, we find that 
\begin{align*}
I_1 (t)= &   \CAL{M}_1(t) \CAL{D}(\zeta _1(t)) \left( \SCR{F} \CAL{M}_2 (t) \SCR{F}^{-1} -1 \right) \hat{w}(t) + \CAL{A}(t), 
\end{align*} 
Hence, we obtain \eqref{1}. 
\subsection{Auxiliary spaces}
Let $(\beta _n, \alpha _n)$ be an admissible, i.e., 
\begin{align} \label{4}
\frac{2}{\beta _n} = n \left( \frac12 - \frac{1}{\alpha _n} \right).
\end{align} 
By defining the notations 
\begin{align*}
\left\| 
F
\right\|_{\alpha , \beta , \gamma , \tau } = \left( \int_{\tau}^{\infty}
\J{s}^{-\gamma} \left\| F(s ,\cdot ) \right\|_{\beta}^{\alpha} ds 
\right) ^{1/ \alpha}, \quad 
\left\| 
F
\right\|_{\infty , \beta , \gamma, \tau } = \sup_{s \geq \tau }
\J{s}^{-\gamma} \left\| F(s ,\cdot ) \right\|_{\beta}. 
\end{align*} 
We define $X_T$ as
\begin{align}  
X_T := \Big\{
\phi \in C\left(
[T, \infty) ; \SCR{S}' ({\bf R} ^n)
\right); \left\|\phi \right\|_{X} &= 
\sup_{\tau \geq T } \tau ^b \| \phi \|_{\infty, 2,\lambda, \tau } \\ & \qquad  + \sup_{ \tau \geq T }  \tau ^{b-2 \lambda} \| \phi \|_{\beta _n, \alpha _n, \lambda, \tau }< \infty
\Big\} \nn 
\end{align}
for $T>r_0$. In particular, we assume that $\lambda >0$ since the case where $\lambda =0$ is very similar to the one considered by \cite{HNW}.


\section{Proof of Theorem \ref{T}}
In this section, we set
\begin{align*}
u(t)-u_p(t)&:=i  \int_{t}^{\infty} U_0(t, s)(F(u(s))-F(u_p(s)))ds+\CAL{E}(t) + \CAL{A}(t),\\
F(\phi(s))&:= \mu |\phi (s) |^{2/(n(1- \lambda))} \phi(s), 
\end{align*}
and consider the following proposition: 
\begin{Prop}\label{P1}
We assume the same assumptions to that in theorem \ref{T}. Then, for $T \geq r_0$, there exists $C>0$ and $\delta_1, \delta_2, \delta _3 >0$ such that 
\begin{align*}
& \left\| u -u_p \right\|_{X_T} 
\\ & \leq  CT^{- \delta _1} \CAL{E}_k(\widehat{u_+}) + C \left\| u-u_p \right\|_{X_T}\left( T^{- \delta_2} \left\| u-u_p \right\|_{X_T}^{\rho_L} + \left\| \widehat{u_+} \right\|_{\infty}^{\rho_L} \right)  +C T^{- \delta _3} \left\|  \widehat{u_+}\right\|_{\infty}^{\rho_L} \| \widehat{u_+} \|_2.
\end{align*}
\end{Prop}

In this section, $\alpha$ denotes a real number such that 
\begin{align*}
2 \alpha < \min(2, 1+ 2/(n(1- \lambda) )-n\lambda /2).
\end{align*}
Here, we remark that for $n=1$, we have $1+2/(1- \lambda) -\lambda /2 > 2 $. As for $n=2$, for all $0< \lambda <1/2$, $1+1/(1- \lambda) -\lambda \geq 2$ holds by $a + 1/a \geq  2$, $a \geq 0$. On the other hand, for $n=3$, 
\begin{align*}
1+ \frac{2}{3(1- \lambda)}  - \frac{3 \lambda }{2}  -2  = \frac{(3 \lambda +1) (3 \lambda -2)}{6(1- \lambda)} \leq 0. 
\end{align*}
holds. This implies that $1+2/(3(1- \lambda))  - 3 \lambda /2  \leq 2$ holds for all $0 < \lambda <1/2$. Consequently, we have 
\begin{align*}
\alpha < 
\begin{cases}
1, & n=1,2, \\ 
1/2 + 1/(n (1- \lambda)) -n \lambda /4, & n=3.
\end{cases}
\end{align*}
Again, we recall that 
\begin{align*}
\lambda (1) < \frac{13-\sqrt{145}}{4}, \quad \lambda (2) < \frac{7-\sqrt{41}}{4} , \quad \lambda (3) < 0.022.
\end{align*}


\begin{Lem}\label{L2}
Let 
\begin{align} \label{8}
0< b <   \lambda + \alpha (1-2\lambda) , \quad  k_1 + 2 \alpha < 1+ 2/(n(1- \lambda)), \quad k_1 = n(1/2-1/\alpha_n).
\end{align}
Then, there exists $0< \delta_1 < \lambda + \alpha (1-2 \lambda) -b $, such that 
\begin{align*}
\left\| R(t) \hat{w}(t) \right\|_{X_T} + \left\| \int_t^{\infty} U_0(t,s) R(s) \hat{w}(s) \frac{ds}{\zeta _2 (s)^{1/(1- \lambda)} } \right\|_{X_T} \leq C T^{- \delta_1} \CAL{E}_{k_1 + 2 \alpha}(\widehat{u_+}).
\end{align*}
\end{Lem}
\Proof{
For simplicity, for the $\tau$-parameter, we use the notation 
\begin{align*}
\| f \|_{\alpha, \beta , \gamma , \tau} = \|  f \|_{\alpha, \beta, \gamma}. 
\end{align*}
We first estimate the term associated with $R(t) \hat{w} (t)$. By straightforward calculations, we have 
\begin{align}\label{2}
\left\| 
R(t) \hat{w}(t) 
\right\|_{\infty,2,\lambda} \leq C \sup_{t \geq \tau} \J{t}^{- \lambda} \left\| (\CAL{M}_2 (t) -1) \SCR{F}^{-1} \hat{w}(t) \right\|_2. 
\end{align}
Together with 
\begin{align*}
\left|  (\CAL{M}_2 (t) -1)  \right| \leq C \left| \frac{\zeta _1 (t)}{\zeta _2 (t)}  \right|^{\alpha} |x|^{2 \alpha} \leq C |t|^{-\alpha (1-2 \lambda)} |x|^{2 \alpha} ,
\end{align*}
we find that $\delta >0$.
\begin{align*}
\left\| 
R(t) \hat{w}(t) 
\right\|_{\infty,2,\lambda} &\leq C \sup_{t \geq \tau} \J{t}^{- \lambda - \alpha (1-2 \lambda)} \left\| \hat{w} (t)\right\|_{H^{2\alpha ,0}} \\ & \leq C \J{\tau}^{-\lambda - \alpha (1- 2 \lambda)} (\log \tau)^{[2\alpha]} 
\CAL{E}_{2 \alpha} ( \widehat{u_+} )
 \\ & \leq C \tau^{-b- \delta} \CAL{E}_{2 \alpha} ( \widehat{u_+} )
 .
\end{align*}
We now estimate $\left\| 
R(t) \hat{w}(t)
\right\|_{\beta _n, \alpha_n,\lambda}$. We remark that for all $f \in L^p({\bf R}^n)$, $p \geq 1$, 
\begin{align*}
\left\| 
\CAL{D}(\zeta _2 (t) ) f 
\right\|_p = \left( 
\int_{{\bf R}^n} \left| 
|\zeta _2 (t)|^{-n/2}f \left( \frac{x}{ \zeta _2 (t)}  \right)
\right| ^p dx
\right)^{1/p} = |\zeta _2 (t)|^{-n(1/2-1/p)} \|f  \|_p
\end{align*}
holds. Here, we employ the Sobolev's embedding, such that 
\begin{align*}
\left\| 
 \hat{w}(t)
\right\|_{\alpha _n} \leq C \left\| \hat{w}(t) \right\|_{H^{k_1,0}} 
\end{align*}
with $k_1> 0$, $2k_1<n$ and 
\begin{align*}
\frac{1}{\alpha _n} = \frac12 - \frac{k_1}{n}.
\end{align*}
Then, by straightforward calculation, 
\begin{align*}
\left\| R(t) \hat{w}(t)\right\|_{\alpha _n} 
&\leq C  |\zeta _2 (t)|^{-n(1/2-1/\alpha _n)}\left\|  \SCR{F} \left( \CAL{M}_2 (t)  -1 \right)\SCR{F}^{-1} \hat{w}(t)  \right\|_{\alpha _n} 
\\ &\leq C \J{t}^{-n(1/2-1/\alpha _n) - (1-2 \lambda) \alpha} (\log t)^{[k_1 + 2 \alpha]}   \left\| \SCR{F}^{-1} \hat{w}(t) \right\|_{0, k_1 + 2 \alpha } \\ & \leq 
C \J{t}^{-n(1/2-1/\alpha _n) - (1-2 \lambda) \alpha} (\log t)^{[k_1 + 2 \alpha]} \CAL{E}_{k_1 + 2 \alpha}( \widehat{u_+}) 
\\ & =: C \gamma (t) \CAL{E}_{k_1 + 2\alpha}( \widehat{u_+})
 .
\end{align*}
Hence, we find 
\begin{align*}
\left\| 
R(t) \hat{w}(t) 
\right\|_{\beta _n, \alpha_n,\lambda} &\leq C  \CAL{E}_{k_1 + 2 \alpha}( \widehat{u_+}) \left( 
\int_{\tau}^{\infty} \J{t}^{- \lambda} ( \gamma (t) )^{\beta _n} dt \right)^{1/\beta _n} . 
\end{align*}
Here, we notice that 
\begin{align*}
& -\lambda - \beta _n \left( n \left(\frac12 - \frac{1}{\alpha _n}  \right) + (1-2\lambda) \alpha   \right) + 1 \\ &= 
-\beta _n \left( 
- (1- \lambda) n \left( \frac14 - \frac{1}{2 \alpha _n} 
\right) + n \left(\frac12 - \frac{1}{\alpha _n}  \right) + (1-2\lambda) \alpha \right) \\ &= 
- \beta _n \left(  \frac{(1+ \lambda)n}{2}  \left( \frac12 - \frac{1}{ \alpha _n} 
\right) + (1-2\lambda) \alpha   \right) <0
\end{align*}
where we use $1/\beta _n = n(1/4 - 1/(2 \alpha _n))$. Since $\alpha _n >2$, we find that 
\begin{align*}
\left\| 
R(t) \hat{w}(t) 
\right\|_{\beta _n, \alpha_n,\lambda}  & \leq 
C \J{\tau}^{- \alpha  (1-2 \lambda  )  -\delta } \CAL{E}_{k_1 + 2 \alpha}( \widehat{u_+}) \\ & \leq C \J{\tau}^{- (b-2 \lambda) - (\lambda + \delta)  } \CAL{E}_{k_1 + 2 \alpha}( \widehat{u_+}) \\ & \leq 
C \J{\tau}^{- (b-2 \lambda) - \delta   } \CAL{E}_{k_1 + 2 \alpha}( \widehat{u_+}) .
\end{align*}

Next, we consider the term associated with $\int U_0(t,s) R(s) \hat{w}(s) s^{-1} ds$. Straightforward calculation and the similar calculation in analyzing the term \eqref{2} shows, 
\begin{align*}
\left\| 
\int_t^{\infty} U_0(t,s) R(s) \hat{w}(s) \frac{ds}{\zeta _2 (s)^{1/(1- \lambda)}}
\right\|_{\infty,2,\lambda}  & \leq C \sup_{t \geq \tau} \J{t}^{- \lambda} \int_t^{\infty} \| R(s) \hat{w}(s) \| \frac{ds}{s} \\ & \leq C  \CAL{E}_{2 \alpha}( \widehat{u_+}) \sup_{t \geq \tau} \J{t}^{- \lambda} \int_{\tau}^{\infty} \J{s}^{ - \alpha(1-2 \lambda)  -1} (\log s)^{[2\alpha]} ds \\ & \leq C \CAL{E}_{2 \alpha}( \widehat{u_+})  \tau^{-b-\delta }  . 
\end{align*}
Moreover, by the Strichartz estimates with an admissible pair $(\infty, 2)$, we have 
\begin{align*}
 \left\| 
\int_t^{\infty} U_0(t,s) R(s) \hat{w}(s) \frac{ds}{\zeta _2(s)^{1/(1- \lambda)} }
\right\|_{\beta _n, \alpha _n,\lambda} & \leq C \left\| R(s) \hat{w}(s) s^{-1} \right\|_{1,2,- \lambda} 
\\ & \leq C  \CAL{E}_{2 \alpha}( \widehat{u_+}) \int_{\tau}^{\infty} \J{s}^{\lambda -1} \J{s}^{- \alpha(1-2 \lambda) } (\log s)^{[2 \alpha]} ds \\ & \leq C \tau^{-(b-2 \lambda)- \delta }  \CAL{E}_{2 \alpha}( \widehat{u_+}), 
\end{align*}
where in order to justify the above, we have to assume that 
\begin{align}\label{3}
\lambda -1 - \alpha (1-2 \lambda) < -1.
\end{align}
For $n=1,2$, we notice that $\alpha < 1$ since $\lambda >0$; hence, \eqref{3} coincides with $\lambda <1/3$. For $n=3$,   
we have 
\begin{align*}
\lambda < \alpha (1 - 2 \lambda) < (1- 2 \lambda) \left( \frac12 + \frac{1}{3(1- \lambda)} - \frac{3}{4} \lambda \right). 
\end{align*}
By straightforward calculation, we have $18 \lambda ^3 -51 \lambda ^2 + 47 \lambda - 10<0$, which will be true for all $0<\lambda < 3/10  \, (>0.02) $.
}


\begin{Lem}\label{lem3.2}
Let $\rho_L:= 2/(n(1- \lambda))$, suppose $1/ \alpha _n < (1- \lambda) /2$ and 
$$  \mathrm{min} (b , \rho_L (b-2 \lambda) )  >1/2 + 2 \lambda - 2\lambda/({\alpha_n} (1 - \lambda)) > 1/2 + \lambda, $$ then there exists $\delta _2 >0$, such that 
\begin{align}\label{str}
&\left\|\int_{t}^{\infty}U_0(t, s)(F(u)-F(u_p))(s)ds\right\|_{X_{T}} 
 \leq C \|u-u_p\|_{X_T} \left( T^{- \delta _2}\|u-u_p\|_{X_{T}}^{\rho_L}+\left\|\widehat{u_+}\right\|_{\infty}^{\rho_L}\right)
\end{align}
holds for any $u \in X_{T}$ with $T\geq r_0$.
\end{Lem}
\Proof{
We employ the same approach proposed by \cite{HNW}. First, we decompose $(F(u)-F(u_p))(s)=F^{(1)}(u(s))+F^{(2)}(u(s))$ with
\begin{align*}
F^{(1)}(u(s))&:=\chi_{\{|u_p(s)|\leq|(u-u_p)(s)|\}}(F(u)-F(u_p))(s),\\
F^{(2)}(u(s))&:=\chi_{\{|u_p(s)|\geq|(u-u_p)(s)|\}}(F(u)-F(u_p))(s). 
\end{align*}
and $\chi_A(x)$ is $1$ or $0$ depending on whether $x\in A$ or not.
Thus, it follows that
\begin{align*}
|F^{(1)}(u)|\leq C|u-u_p|^{\rho_L}|u-u_p|,\quad |F^{(2)}(u)| \leq C|u_p|^{\rho_L}|u-u_p| .
\end{align*}
We estimate 
\begin{align*}
|(F(u(s))-F(u_p(s)))|&\leq |F^{(1)}(u(s))|+|F^{(2)}(u(s)) 
\\ & \leq C|(u-u_p)(s)|^{\rho_L}|(u-u_p)(s)|  +C|u_p(s)|^{\rho_L}|(u-u_p)(s)|.
\end{align*}
Hereafter, we estimate \eqref{str} with respect to $F^{(1)}(u(s))$ and $F^{(2)}(u(s))$. \\ {\bf Estimation for $F^{(1)}$.} \\ 
By the Strichartz estimate, the term associated with $F^{(1)}(u(s))$ is estimated as 
\begin{align*}
\left\|\int_t^{\infty}U_0(t, s)F^{(1)}(u(s))ds\right\|_{\infty, 2,\lambda}
&\leq C \left\||u-u_p|^{1+\rho_L}\right\|_{\rho ', q', -\lambda}\\
&\leq C \left(\int_{\tau}^{\infty}\langle s \rangle^{\lambda}\left\||(u-u_p)(s)|^{1+\rho_L}\right\|_{q'}^{\rho'}ds \right)^{1/\rho'}, 
\end{align*}
where $\tau \geq T$. Since it holds that
\begin{align*}
\left\||(u-u_p)(s)|^{1+\rho_L}\right\|_{q'}&\leq \left\|(u-u_p)(s)\right\|_{(2q' \rho_L)/(2-q')}^{\rho_L}\times \|(u-u_p)(s)\|_{2},
\end{align*}
by Holder's inequality with $\alpha_n=(2q'\rho_L)/(2-q')$, we have
\begin{align}
& \nn \left(\int_{\tau}^{\infty}\langle s \rangle^{\lambda}\left\||(u-u_p)(s)|^{1+\rho_L}\right\|_{q'}^{\rho'}ds \right)^{1/\rho'}\\
& \nn \leq C\left(\int_{\tau}^{\infty}\langle s \rangle^{\lambda}\left\|(u-u_p)(s)\right\|_{\alpha_n}^{\rho'\rho_L}\times \|(u-u_p)(s)\|_{2}^{\rho'}d s \right)^{1/\rho'}\\
& \nn \leq C\left\|\langle s \rangle^{(1+1/\rho')\lambda -b 
}\left\|(u-u_p)(s)\right\|_{\alpha_n}^{\rho_L}\right\|_{L_{\tau}^{\rho'}}\left(\sup_{s\in[\tau ,\infty)} \langle s \rangle^{b} \langle s \rangle^{-\lambda}\|(u-u_p)(s)\|_{2}\right)\\
& \nn \leq C\left\|\langle s \rangle^{(1+1/\rho'+\rho_L/\beta_n)\lambda-b}\right\|_{L^p_{\tau}} 
\left(\int_t^{\infty}\langle s \rangle^{-\lambda}\left\|(u-u_p)(s)\right\|_{\alpha_n}^{\beta_n}ds \right)^{\rho_L/\beta_n} \\
& \label{6} \qquad \times \sup_{\tau \geq T } \left(  \tau^{b}\|u-u_p\|_{\infty, 2,\lambda} \right), 
\end{align}
where $p>1$ satisfies $\rho_L/\beta_n+1/p=1/\rho'$, 
$$\| f(s) \|_{L^{\theta}_\tau} := \left( \int_{\tau}^{\infty} |f(s)|^{\theta} ds \right)^{1/ \theta}, $$ 
and we use 
\begin{align*}
\sup_{s\in[\tau ,\infty)} \langle s \rangle^{b} \langle s \rangle^{-\lambda}\|(u-u_p)(s)\|_{2} 
& \leq \sup_{s\in[\tau ,\infty)} \langle s \rangle^{b} \sup_{\kappa \geq s} \langle \kappa  \rangle^{-\lambda}\|(u-u_p)(\kappa)\|_{2} \\ & \leq \sup_{s\in[\tau ,\infty)} \langle s \rangle^{b} \|u-u_p\|_{\infty, 2, \lambda,s}  \\  & \leq C  \sup_{(s \geq) \tau \geq T } \left(  \tau^{b}\|u-u_p\|_{\infty, 2,\lambda} \right).
\end{align*} 
Because of the condition $n(1/ q' -1/2) =2/ \rho < 1$, the restriction 
\begin{align*}
\alpha _n > n \rho_L = \frac{2}{1- \lambda}
\end{align*}
appears. Then, taking $b>0$, such that
\begin{align} \label{7}
b>1/2 + 2 \lambda - 2\lambda/({\alpha_n} (1 - \lambda)) > 1/2 + \lambda
\end{align}
holds.
\begin{align*}
\left\|\langle s \rangle^{(1+1/\rho'+\rho_L/\beta_n)\lambda-b}\right\|_{L^p_{\tau}}&\leq C\left(\int_{\tau}^{\infty} s^{((1+1/\rho'+\rho_L/\beta_n)\lambda-b)p}d s\right)^{1/p}\\
&\leq C \tau^{1/p+(1+1/\rho'+\rho_L/\beta_n)\lambda-b}\\
&= C \tau^{1/2 + 2 \lambda - 2\lambda/({\alpha_n} (1 - \lambda))-b} 
\end{align*}
is bounded, where we use \eqref{4} and 
\begin{align*}
\frac{1}{\rho '} = 1- \frac{n}{2} \frac{2-q'}{2q'} = 1- \frac{n \rho_L}{2 \alpha _n}.
\end{align*}
Therefore, we obtain that \eqref{6} is smaller than 
\begin{align*}
&\left\|\langle s \rangle^{(1+1/\rho'+\rho_L/\beta_n)\lambda-b}\right\|_{L^p_{\tau}} \left\| u-u_p \right\|_{\beta _n, \alpha _n, \lambda}^{\rho_L} 
\sup_{\tau\in[T,\infty)}\tau ^{b}\|u-u_p\|_{\infty, 2,\lambda}\\
&\leq C\tau ^{1/2 + 2 \lambda - 2\lambda/({\alpha_n} (1 - \lambda))-b} 
\|u-u_p\|_{\beta _n, \alpha _n,\lambda}^{\rho_L}
\|u-u_p\|_{X_T}\\
&\leq C \tau ^{1/2 + 2 \lambda - 2\lambda/({\alpha_n} (1 - \lambda))-b- \rho_L(b-2 \lambda)} 
\left( \tau ^{b-2\lambda}\|u-u_p\|_{\beta _n, \alpha _n,\lambda}\right)^{\rho_L}
\|u-u_p\|_{X_T}\\
&\leq C\tau ^{-b}\|u-u_p\|_{X_T}^{\rho_L +1 } \times  \tau ^{1/2 + 2 \lambda - 2\lambda/({\alpha_n} (1 - \lambda))-\rho_L (b-2\lambda )}
.
\end{align*}
We suppose that 
\begin{align}\label{9}
\rho_L (b -2\lambda )  >1/2 + 2 \lambda - 2\lambda/({\alpha_n} (1 - \lambda)) > 1/2 + \lambda. 
\end{align}
We find, 
\begin{align*}
\sup_{\tau \geq T } \tau ^{b}\left\|\int_t^{\infty}U_0(t, s)F^{(1)}(u(s))ds\right\|_{\infty, 2,\lambda}
 \leq C  T ^{1/2 + 2 \lambda - 2\lambda/({\alpha_n} (1 - \lambda))-\rho_L (b-2\lambda )}\left\| u- u_p \right\|_{X_T}^{\rho_L +1}, 
\end{align*}
where we remark that in \eqref{7}, $\lambda < 1/3$ and $ 2/(1- \lambda) < \alpha _n < \infty  $, 
\begin{align*} 
b-2 \lambda > 1/2 + \lambda >0 .
\end{align*} By the same calculation, we find there exists $\delta _2 >0$, such that
\begin{align*}
&\sup_{\tau \geq T } \tau^{b-2 \lambda} \left\| 
\int_t^{\infty } U_0(t,s) F^{(1)} (u(s)) ds 
\right\|_{\beta _n, \alpha _n, \lambda}\\ &  \leq C \left\| u-u_p \right\|_{X_T}^{\rho_L +1} \times \sup_{\tau \geq T }  \tau^{1/2 + 2 \lambda - 2\lambda/({\alpha_n} (1 - \lambda))-(\rho_L -1) (b-2\lambda) -b} \\ & \leq  C T^{-\delta _2}\left\| u-u_p \right\|_{X_T}^{\rho_L +1}, 
\end{align*} 
where we use $(b-2 \lambda) -b <0$. \\ ~~ \\ 
{\bf Estimation for $F^{(2)}$.} \\ 
First, remarking that $2/(1- \lambda)< \alpha _n < \infty$,
\begin{align*}
-b < - 1/2 -\lambda < - \lambda , \quad \mbox{and} \quad -b+2 \lambda <0.
\end{align*}
Then
\begin{align*}
& \left\| 
\int_t^{\infty} U_0(t,s) F^{(2)} (u(s)) ds 
\right\|_{\infty,2,\lambda} \\ & = \sup_{t  \geq \tau } \J{t}^{- \lambda} \left\| \int_{t}^{\infty} U_0(t,s) F^{(2)} (u(s)) ds \right\|_2 \\ & \leq C \sup_{ t \geq \tau} \J{t}^{- \lambda} \int_t^{\infty} \left\| u_p (s) \right\|_{\infty}^{\rho_L} \left\| u(s) -u_p(s) \right\|_2 ds \\ & \leq  C \|  \widehat{u_+} \|_{\infty}^{\rho_L}\sup_{t \geq \tau } \J{t}^{- \lambda}  \sup_{s \geq t} \left( \J{s}^{-\lambda +b} \| u(s) -u_p(s) \|_2 \right)  \int_t^{\infty} \J{s}^{-1-b + \lambda} ds \\ & \leq C \tau ^{-b} \|  \widehat{u_+} \|_{\infty}^{\rho_L} \left\| u-u_p \right\|_{X_T} 
,
\end{align*}
where we use 
\begin{align*}
\| u_p(s) \|_{\infty}^{\rho_L} \leq C |\zeta _2 (s)|^{\rho_L} \|  \widehat{u_+} \|_{\infty} \leq C s^{-1}  \|  \widehat{u_+} \|_{\infty}. 
\end{align*} 
Strichartz estimates with an admissible pair $(\infty,2)$ yield 
\begin{align*}
\left\|\int_t^{\infty}U_0(t, s)F^{(2)}(u(s))ds\right\|_{\beta _n, \alpha _n,\lambda}
& \leq C \left\|  | {u_p}|^{\rho_L} | (u-u_p) | \right\|_{1,2, -\lambda} \\ 
&\leq C  \int_{\tau}^{\infty}\langle s \rangle^{\lambda}\left\| u_p (s)\right\|^{\rho_L}_{\infty} \| (u-u_p)(s) \|_2 d s .
\end{align*}
By the same argument as above, we find that the last term of the above inequality is smaller than 
\begin{align*}
&C \| u-u_p \|_{X_T} \| \widehat{ u_+} \|_{\infty}^{\rho_L} \int_{\tau}^{\infty} \J{s}^{-b + 2 \lambda -1} d s \\ & \leq C \tau^{-b + 2 \lambda} \| u-u_p \|_{X_T} \left\| \widehat{ u_+} \right\|_{\infty}^{\rho_L}.
\end{align*}

}

\begin{Lem}\label{L-ad1}
Under the assumption \ref{A1} and for $b < \lambda + \alpha (1- 2 \lambda)$, there exists $\delta _3 >0$ such that 
\begin{align*}
\left\| \CAL{A} (t) \right\|_{X_T} \leq C T^{-\delta _3} \| \widehat{u_+} \|_{\infty}^{\rho_L} \| \widehat{u}_+ \|_2
\end{align*}
\end{Lem}
\Proof{
Since $\zeta _2 (s) = c_{2,+} s^{1- \lambda} + \CAL{O}(t^{ \lambda})$, we have $s \leq  (\zeta _s (s) /c_{2,+})^{1/(1- \lambda)} + \CAL{O}(t^{ \lambda /(1- \lambda)}) $, and which implies 
\begin{align*}
\left|  
\frac{c_+ s}{\zeta _2 (s)^{1/(1- \lambda)}}  -1
\right| \leq C s^{\lambda /(1- \lambda) -1}, 
\end{align*}
where we use $c_+ = c_{2,+}^{1/ (1- \lambda)}$. Then 
\begin{align*}
\left\|\CAL{A}(t) \right\|_{\infty, 2, \lambda} &\leq C  \| \widehat{u_+} \|_{\infty}^{\rho_L} \| \widehat{u_+} \|_2 \sup_{t\geq \tau }  \J{t}^{- \lambda} \int_{t}^{\infty} \J{s}^{-2  + \lambda/(1- \lambda)} ds \\ & \leq C  \| \widehat{u_+} \|_{\infty}^{\rho_L} \| \widehat{u_+} \|_2 \tau^{-b + ( b -1 - \lambda   + \lambda/(1- \lambda) )}
\end{align*}
holds. Here we remark that 
\begin{align*}
b -1-\lambda + \frac{\lambda}{1- \lambda} < \alpha (1- 2 \lambda) -1 + \frac{\lambda}{1- \lambda} <  \frac{\lambda}{1- \lambda}-2 \lambda =- \frac{\lambda (1- \lambda)}{1- \lambda} <0. 
\end{align*}
On the other hand, 
\begin{align*}
\left\|\CAL{A}(t) \right\|_{\beta _n, \alpha _n , \lambda} & = \left\|\int_t^{\infty} U_0 (t,s) U_0 (s,0) \SCR{F}^{-1} \left( 
\frac{c_+ s}{\zeta _2 (s)^{1/(1- \lambda)}} -1
\right) F(\hat{w}(s)) \frac{ds}{c_+ s}  \right\|_{\beta _n, \alpha _n, \lambda} \\ 
& \leq C \int_{\tau}^{\infty} \J{s}^{\lambda} \left\| U_0 (s,0) \SCR{F}^{-1} \left( 
\frac{c_+ s}{\zeta _2 (s)^{1/(1- \lambda)}} -1
\right) F(\hat{w}(s))  \right\|_{2} \frac{ds}{ s} 
\\
&\leq C  \| \widehat{u_+} \|_{\infty}^{\rho_L} \| \widehat{u_+} \|_2 \int_{\tau}^{\infty}  \J{s}^{ -2 + \lambda + \lambda /(1- \lambda) } ds
\\
&\leq C  \| \widehat{u_+} \|_{\infty}^{\rho_L} \| \widehat{u_+} \|_2  \J{\tau}^{ -(b -2\lambda) + (b -1 - \lambda + \lambda /(1- \lambda)) } .
\end{align*}
These prove Lemma \ref{L-ad1}.
}


\subsection{Proof of Theorem \ref{T}}

By Lemma \ref{L2} with $1/\alpha_n < (1- \lambda) /2$, i.e., $k_1 > n\lambda /2$, lemma \ref{lem3.2}, lemma \ref{L-ad1} and proposition \ref{P1}. By considering sufficiently small $\| \widehat{u_+} \|_{\infty}$ and sufficiently large $T >0$ compared with the constant $C>0$, there exists $0< \ep _1 \ll 1$ and $\delta _4 >0$, such that 
\begin{align*}
\left\|  u-u_p \right\|_{X_T} \leq \ep_1 + \left\|  u-u_p \right\|_{X_T}^{1+ \delta _4} .
\end{align*}
The bootstrap argument yields 
\begin{align*}
\left\|  u-u_p \right\|_{X_T} \leq (1 + \delta _4) \ep_1/\delta _4, 
\end{align*}
which implies theorem \ref{T}. 

\subsection{Existence of an admissible pair $(\beta _n, \alpha _n)$}
Here, we show the existence of an admissible pair $(\beta _n, \alpha _n)$ that satisfies conditions \eqref{8}, \eqref{7}, and \eqref{9}. It suffices to find the condition of $\lambda$, such that
\begin{align}\label{10}
\frac{1}{2} +  \lambda < b < \lambda + \alpha (1-2\lambda) 
\end{align}
and 
\begin{align} \label{12}
2\lambda + \frac{1}{2 \rho_L} + \frac{ \lambda}{\rho_L}  < b < \lambda + \alpha (1-2\lambda) .
\end{align}
For $n=1,2$, we have $\alpha <1$, and thus, \eqref{10} implies that $\lambda (n) < 1/4$. On the other hand, for $n=3$, noting $\alpha < 1/2 + 1/(3(1- \lambda)) - 3\lambda /4 $, \eqref{10} implies
\begin{align*}
18 \lambda ^3 -39 \lambda ^2 + 29 \lambda -4 < 0 \Rightarrow \lambda (3) < 0.176...
\end{align*}

Next, we consider the condition of $\lambda$ such that \eqref{12} will be true; for $n=1$, then $\rho_L = 2/(1- \lambda)$ and $\alpha <1$ yield 
\begin{align*}
2 \lambda ^2 -13 \lambda + 3 > 0 \Rightarrow \lambda (1) < \frac{13 - \sqrt{145}}{4} = 0.239... < \frac14.
\end{align*}
For $n=2$, $\rho_L = 1/(1- \lambda)$ and $\alpha <1$ yield 
\begin{align*}
2 \lambda ^2 -7 \lambda + 1 >0 \Rightarrow \lambda (2) < \frac{7-\sqrt{41}}{4} = 0.149....
\end{align*}
For $n=3$, $\rho_L = 2/(3(1- \lambda))$, and $\alpha < 1/2 + 1/(3(1- \lambda))-3 \lambda /4$ yield 
\begin{align*}
36 \lambda ^3 -78 \lambda ^2 + 47 \lambda -1 <0 \Rightarrow \lambda (3) < 0.022...
\end{align*}

~~ \\ ~~ \\ 
{\bf Acknowledgement:} \\ 
The author is partially supported by the Grant-in-Aid for Young Scientists \#20K14328 from JSPS. The author is grateful to Professor Kota Uriya and Ryo Muramatsu for giving some valuable comments.



\begin{thebibliography}{ABCD}

\bibitem{Ca} Carles, R.: Nonlinear Schr\"{o}dinger equation with time dependent potential, Comm. Math. Sci., \textbf{9}, 937--964 (2011).







\bibitem{GO} J. Ginibre, T. Ozawa, Long range scattering for nonlinear Schr\"{o}dinger and Hartree equations in space dimension $n \geq 2$, Com. Math. Phys., \textbf{151}, 619--645 (1993). 


\bibitem{HN} N. Hayashi, P. I. Naumkin, Asymptotics for large time of solutions to the nonlinear Schr\"{o}dinger and Hartree equations, American J. of Math., \textbf{120}, 369--389 (1998). 

\bibitem{HNW} N. Hayashi, P. I. Naumkin, H. Wang, Modified wave operators for nonlinear Schr\"{o}dinger equations in lower order Sobolev spaces, J. Hyperbolic Diff. Eqn., \textbf{8}, 759--775 (2011)

\bibitem{HO} N. Hayashi, T. Ozawa, Scattering theory in the weighted $L^2(\mathbb{R}^n)$ spaces for some Schr\"{o}dinger equations, Ann. Inst. H. Poincar\'{e}, Phys. Th\'{e}or., \textbf{48}, 17--37 (1988)

\bibitem{IK} A. Ishida, M. Kawamoto, Existence and nonexistence of wave operators for time-decaying harmonic oscillators, Rep. Math. Phys., \textbf{85}, 335--350 (2020). 

\bibitem{Ka} M. Kawamoto, Quantum scattering for time-decaying harmonic oscillators, arXiv:1704.03714.

\bibitem{Ka2} M. Kawamoto, Mourre theory for time-periodic magnetic fields, J. Funct. Anal., \textbf{277}, 1--30 (2019).  

\bibitem{Ka3} M. Kawamoto, Strichartz estimates for Schr\"{o}dinger operators with square potential with time-dependent coefficients, To appear in Diff. Eqn. Dyn. Sys. (arXiv:1805.07991)

\bibitem{KY} M. Kawamoto, T. Yoneyama, Strichartz estimates for harmonic potential with time-decaying coefficient, J. Evol. Eqn. \textbf{18}, 127--142 (2017) .

\bibitem{KM} M. Kawamoto, R. Muramatsu, Asymptotic behavior of solutions to nonlinear Schr\"{o}dinger equations with time-dependent harmonic potentials, J. Evol. Eqn. \textbf{21}, 699--723, (2021).

\bibitem{Ko} L. E. Korotyaev, On scattering in an external,
	    homogeneous, time-periodic magnetic field
	    Math. USSR-Sb., \textbf{66}, 499--522 (1990). 
	    
	    
\bibitem{MM} S. Masaki, H. Miyazaki, Long range scattering for nonlinear Schr\"{o}dinger equations with critical homogeneous nonlinearity, SIAM, J. Math. Anal. \textbf{50}, 3251--3270 (2018). 

\bibitem{MMU} S. Masaki, H. Miyazaki, K. Uriya, Long-range scattering for nonlinear Schr\"{o}dinger equations with critical homogeneous nonlinearity in three space dimensions, Transaction of the A.M.S., \textbf{371}, 7925--7947 (2019).    
	    
	    
\bibitem{O} T. Ozawa, Long range scattering for nonlinear Schr\"{o}dinger equations in one space dimension, Com. Math. Phys., \textbf{139}, 479--493 (1991). 





\bibitem{V} M. Visan, The defocusing energy-critical nonlinear Schr\"{o}dinger equation in higher dimensions, Duke Math. J., \textbf{138}, 281--374 (2007).







\end{thebibliography}
\end{document}